\documentstyle[12pt,a4]{article}
\begin{document}
\thispagestyle{empty}
\begin{center}
{\large\bf THE NONCOMMUTATIVE INHOMOGENEOUS HOPF ALGEBRA}\\
\vskip2truecmM. Lagraa and N. Touhami\\
\vskip1truecm Laboratoire de Physique Th\'eorique\\Universit\'e d'Oran Es-S\'enia, 31100, Alg\'erie\\
\vskip1truecm
\end{center}
\vskip3truecm
\begin{abstract}
From the bicovariant first order differential calculus on 
inhomogeneous Hopf algebra ${\cal B}$ we construct the set of right-invariant 
Maurer-Cartan  one-forms considered as a right-invariant basis of a bicovariant ${\cal B}$-bimodule over which we develop the Woronowicz's general theory of
differential calculus on quantum groups. In this formalism, we introduce 
suitable functionals on ${\cal B}$ which control the inhomogeneous commutation 
rules. In particular we find that the homogeneous part of commutation rules 
between the translations and those between the generators of the homogeneous
part of ${\cal B}$ and translations are controled by different $R$-m
atrices satisfying nontrivial characteristic equations.
\end{abstract}
\newpage
\section{Introduction}
The Poincare group plays a fundamental role in physics. It is intrinsically 
connected to the geometry of the space-time on which the
 physical systems are described and could be invariant under its action. Then,  it is especially interesting to study the non-commutative version of the 
Poincare group from which on can hope to obtain new insight on the underling 
space-time geometry as, for example, an improved ultraviolet in quantum field 
theories or a description of symmetries in a future quantum Einstein-Cartan 
gravity.\\ 
The construction of quantum Poincare group and quantum space-times has alreadybeen considered by several authors [1]-[8]. These constructions 
starteither fromthe existence of a $R$-matrix and the consistency of the 
commutation rules betweenthe different elements of the generators of the 
Poincare group and those of the quantum space-times or from the projection 
of homogeneous quantum groups. In certain works this consistency is only 
obtained by introducing extra generator [4], [6].\\
Recently, in the work P. Podles and S.L. Woronowicz [9] inhomogeneous 
commutation rules for inhomogeneous Hopf algebra without dilatation have 
been constructed. 
This construction is based on the existence of a bicovariant subimodule of the 
inhomogeneous algebra regarded as a left module (over the homogeneous part) 
generated by translations and leads to a same $R$-matrix, subject to the 
condition $R^{2}=I^{\otimes2}$, which controls the homogeneous parts of the 
different commutation rules.\\ In this paper we present a quite different 
method based on the differential calculus on the inhomogeneous Hopf algebra 
from which we derive the corresponding inhomogeneous commutation rules. 
In this framework, we show that the homogeneous parts of the commutation rules 
are controled by different $R$-matrices satisfying nontrivial characteristic 
equations. Although our formalism is quit in [9], it gives, under special 
conditions, identical results.\\The present paper is organized as follows: 
In section 2, we recall some basic notions about the inhomogeneous Hopf algebra and construct a right-invariant basis of the bimodule over this algebra. 
This right-invariant basis allows us to generalize, In section 3, the 
Woronowicz proposal [10] to the inhomogeneous Hopf algebra by introducing 
suitable functionals over ${\cal B}$. The study ofthese functionals leads 
us to construct inhomogeneous commutation rules between the elements of 
${\cal B}$ whose homogeneous parts are controled by different 
$R$-matrices satisfying nontrivial characteristic equations. Finally, 
we investigate in section 4 the consistency conditions between these 
functionals and the different commutation rules of the inhomogeneous 
quantum group.\\Throughout this paper we use Einstein's convention 
(sum over repeated indices) and for the simplicity of calculations we define 
$v \top w \in M_{NN}({\cal B})$ by 
$(v \top w)^{nm}_{~~k\ell} = v^{n}_{~k}w^{m}_{~\ell}$, 
$n,m,k, \ell \in 1, \cdots ,N$, for any $v$ and $w \in M_{N}({\cal B})$.\\
\section{Differential Calculus On Inhomogeneous \\ Quantum Groups}
In this section, we start by recalling some basic background about 
inhomogeneousHopf algebra ${\cal B}$ and the covariant first order 
differential calculus to  construct a bicovariant bimodule $\Gamma$ 
of one-forms over ${\cal B}$ and a    basis of the vector space 
$\Gamma_{inv} \subset \Gamma$ of all right invariant  elements of 
$\Gamma$.\\ 
An inhomogeneous quantum group ${\cal G}$ is built from a quantum 
group${\cal H}$ and translations described by elements $p^{n},n=1,...,N$ 
corresponding to an irreducible representation $\Lambda^{n}_{~m}$ of ${\cal H}$. The corresponding Hopf algebras are treated as algebras of functions on quantum groups, Poly$({\cal G})$ = ${\cal B}$ and Poly$({\cal H})$ = ${\cal A}$. 
More precisely, following [9], one defines ${\cal B}$ as:\\ 1. An abstract
unital H
${\cal A}$ and  elements $p^{n}$ such that $I_{\cal B} = I_{\cal A} = I $.\\ 
2. ${\cal A}$ is a Hopf sub-algebra of ${\cal B}$.\\
3. ${\cal P} = \left( \begin{array}{cc} \Lambda  &  p \\ 0  &  I \end{array}    \right)$
is a representation of ${\cal G}$.\\4. There exists $n \in 1, \cdots ,N$ such 
that $p^{n} \not\in {\cal A}$.\\5. $\Gamma_{p} {\cal A} \subset \Gamma_{p}$
 where $\Gamma_{p} = {\cal A}X+{\cal A}$,  $X=span(p^{n}, n=1, \cdots ,N)$.\\
By virtue of 2.-3., ${\cal B}$ is endowed with the following linear maps:\\
- The coaction (algebra homomorphism) $\Delta : {\cal B}\rightarrow{\cal B} 
\otimes {\cal B}$ \\
\begin{eqnarray*}
\Delta(a) = a_{(1)} \otimes a_{(2)},    &a \in {\cal B}
\end{eqnarray*}
satisfying the coassociativity condition\\
\begin{eqnarray*}
(\Delta \otimes id)\Delta(a) = (id \otimes \Delta)\Delta(a) = a_{(1)} \otimes 
a_{(2)} \otimes a_{(3)}.
\end{eqnarray*}
Due to the 
representation ${\cal P}$ of ${\cal G}$ the coaction acts on the generators as\\\begin{eqnarray*}
\Delta (\Lambda^{n}_{m})&=&\Lambda^{n}_{~k} \otimes \Lambda^{k}_{~m}, \\ 
\Delta(p^{n})&=&\Lambda^{n}_{~k} \otimes p^{k} + p^{n} \otimes I.
\end{eqnarray*}
- The counit (character) $\varepsilon : {\cal B} \rightarrow {\cal C}$ 
satisfying\\
\begin{eqnarray*}
(\varepsilon \otimes id)\Delta(a) = (id \otimes \varepsilon)\Delta(a) =a, 
&a \in {\cal B},
\end{eqnarray*}
acts on the generators as\\
\begin{eqnarray*}
\varepsilon( \Lambda^{n}_{~m}) = \delta^{n}_{m}  &and~~ \varepsilon (p^{n}) = 0 .\end{eqnarray*}
- The antipode (algebra antihomomorphism) $S:{\cal B} \rightarrow {\cal B}$
satisfying\\
\begin{eqnarray*}
m \circ (S \otimes id)\Delta(a) = m \circ (id \otimes S)\Delta(a) =
I_{\cal B}\varepsilon(a)  &,~~a\in {\cal B},
\end{eqnarray*}
where $m:{\cal B} \otimes {\cal  B} \rightarrow {\cal B}$ is the multiplication map. it acts on the generators as\\
\begin{eqnarray*}
S(\Lambda^{n}_{k})\Lambda^{k}_{m}&=&\varepsilon(a)^{n}_{m}= \Lambda^{n}_{k}S(\Lambda^{k}_{m}) , \\S(\Lambda^{n}_{k})p^{k}+S(p^{n})&=&\varepsilon (p^{n}) =0= \Lambda^{n}_{k}S(p^{k})+p^{n}.
\end{eqnarray*}
One says [10] that (${\Gamma},d$) is a first order differentialcalculus on the inhomogeneous Hopf algebra (${\cal B},\Delta, S, \varepsilon)$if $d:{\cal B}\rightarrow \Gamma$ is a linear map obeying the Leibniz rule, $d(ab)=(da)b+a(db)$ for any $a,b \in {\cal B}$. $\Gamma$ is a bimodule over ${\cal B}$ and every element of $\Gamma$ is of the form $\sum_{k} a_{k}db_{k}$,where $a_{k}, b_{k} \in {\cal B}$.We say that (${\Gamma},d$) is left-covariant if for any $a_{k},b_{k} \in {\cal B}$\\
\begin{eqnarray*}
\sum_{k} a_{k}db_{k} = 0 \Leftrightarrow \sum_{k} \Delta (a_{k})(id \otimes d)\Delta(b_{k})=0
\end{eqnarray*}
and is right-covariant if
\begin{eqnarray*}
\sum_{k} a_{k}db_{k} = 0 \Leftrightarrow \sum_{k} \Delta (a_{k}) (d \otimes id) \Delta (b_{k}) = 0.
\end{eqnarray*}
(${\Gamma},d$) is bicovariant if it is left- an
d right-covariant. This notionof covariant differential calculus leads to the left-coaction$\Delta_{L}$ and right-coaction $\Delta_{R}$ which are bimodule homomorphisms\\
\begin{eqnarray}
\Delta_{L(R)} (a\rho b)= \Delta(a) \Delta_{L(R)} (\rho ) \Delta (b),    &a,b \in {\cal B},\rho \in \Gamma
\end{eqnarray}
and satisfy\\
\begin{eqnarray}
\Delta_{L} d=(id \otimes d) \Delta & ,~~ \Delta_{R} d=(d \otimes id) \Delta ,\nonumber\\( \varepsilon \otimes id) \Delta_{L} = id &and~~(id \otimes \varepsilon ) \Delta_{R}=id.
\end{eqnarray}
For a bicovariant bimodule, the following condition is satisfied
\begin{eqnarray*}
( \Delta_{L} \otimes id) \Delta_{R} = (id \otimes \Delta_{R} ) \Delta_{L}.
\end{eqnarray*}
Since, ${\cal B} = Poly({\cal G}) = Poly( \Lambda , p)$, $d \Lambda^{n}_{~m}$and $dp^{n}$ generate $\Gamma$ as a bimodule of one-forms over ${\cal B}$, we have\\{\bf Lemma (2,1)}\\
\begin{eqnarray}
\Theta^{n}_{Rm} = d \Lambda^{n}_{~k} S( \Lambda^{k}_{m})  &and~~\Pi^{n}_{R} = dp^{n} - \Theta^{n}_{~k}p^{k}
\end{eqnarray}
the vector space $\Gamma_{inv} \subset \Gamma$ of all right-invariant elements of $\Gamma$.\\\vskip1truecm
$proof.$ Since $p^{n}$ correspond to an irreducible representation $ \Lambda $ of ${\cal H}$, they are linearly independent, so are  $S(p^{n})$ = $- S(\Lambda^{n}_{~m})p^{m}$. Then $d \Lambda^{n}_{~m}$ and $dp^{n}$or $d \Lambda^{n}_{~m}$ and $dS(p^{n})$ form a basis of the one-form vector space which generates $ \Gamma$ as a ${\cal B}$-bimodule from which we can choose an another basis $d( \Lambda^{n}_{~k})S( \Lambda^{k}_{~m}) = \Theta^{n}_{Rm}$ and $ \Lambda^{n}_{~k}dS(p^{k}) = p^{n} - \Theta^{n}_{Rm}p^{m} = \Pi^{n}_{R}$. To show the rightinvariance, it suffices to apply (1), (2) and the properties of the Hopf algebra structure of ${\cal B}$ to get \\
\begin{eqnarray}
\Delta_{R}( \Theta^{n}_{Rm})&=&d \Lambda^{n}_{~\ell}S( \Lambda^{p}_{~m}) \otimes \Lambda^{\ell}_{~k}S( \Lambda^{k}_{~p}) = \Theta^{n}_{Rm} \otimes I ,\\
\Delta_{R}( \Pi^{n}_{R})&=&d \Lambda^{n}_{~k} \otimes p^{k} + dp^{n} \otimes I - \Theta^{n}_{Rm} \Lambda^{m}_{~k} \otimes p^{k} - \Theta^{n}_{Rm}p^{m} \otimes I \nonumber \\
 &=&\Theta^{n}_{Rm} \Lambda^{m}_{~k} \otimes p^{k} + (dp^{n} - \Theta^{n}_{Rm}p^{m}) \otimes I - \Theta^{n}_{Rm} \Lambda^{m}_{~k} \otimes p^{k} = \Pi^{n}_{R} \otimes I ,\nonumber \\ \\
\Delta_{L}( \Theta^{n}_{Rm})&=&\Lambda^{n}_{~\ell} S( \Lambda^{k}_{~m}) \otimes \Theta^{\ell}_{Rk} ,\\
\Delta_{L}( \Pi^{n}_{R}) &=& \Lambda^{n}_{~k} \otimes dp^{k} - \Lambda^{n}_{~k} S( \Lambda^{\ell}_{~m}) \Lambda^{m}_{~p} \otimes \Theta^{k}_{R\ell} p^{p} - \Lambda^{n}_{~k} S( \Lambda^{\ell}_{~m})p^{m} \otimes \Theta^{k}_{R\ell} \nonumber \\&=&\Lambda^{n}_{~k} \otimes (dp^{k} - \Theta^{k}_{R\ell} p^{\ell}) +\Lambda^{n}_{~k} S(p^{\ell}) \otimes \Theta^{k}_{R\ell} \nonumber \\
&=&\Lambda^{n}_{~k} \otimes \Pi^{k}_{R} + \Lambda^{n}_{~k} S(p^{\ell}) \otimes \Theta^{k}_{R\ell}.
\end{eqnarray}
Q.E.D.\\
$Remark (2,1)$: Note that we can construct from $\Pi^{n}_{R}$ and $\Theta^{n}_{Rm}$ a left-invariant basis as:
\begin{eqnarray}
_{Lm}&=&S( \Lambda^{n}_{~k} S( \Lambda^{\ell}_{~m})) \Theta^{k}_{R\ell} ,\\
\Pi^{n}_{L}&=&S( \Lambda^{n}_{~m}) \Pi^{m}_{R} + S( \Lambda^{n}_{~m} S(p^{k})) \Theta^{m}_{Rk}
\end{eqnarray}
satisfying\\
\begin{eqnarray*}
\Delta_{L} ( \Theta^{n}_{Lm})&=&I \otimes \Theta^{n}_{Lm}~~and~~  \Delta_{L} ( \Pi^{n}_{L}) = I \otimes \Pi^{n}_{L} ,\\\Delta_{R} ( \Theta^{n}_{Lm})&=&\Theta^{k}_{L\ell} \otimes S( \Lambda^{n}_{~k}S( \Lambda^{\ell}_{~m})) ,\\\Delta_{R} ( \Pi^{n}_{L})&=&\Pi^{k}_{L} \otimes S( \Lambda^{n}_{~k}) + \Theta^{k}_{L\ell} \otimes S( \Lambda^{n}_{~k} S(p^{l})).
\end{eqnarray*}
The bicovariance condition of these basis can be checked by direct computation. In the following, we consider the basis $\Theta^{n}_{Rm} = \Theta^{n}_{~m}$ and $\Pi^{n}_{R} = \Pi^{n}$. 
\section{\bf Commutation Rules For Inhomogeneous Hopf Algebras}Now, we are ready to study the commutation rules for inhomogeneous Hopf algebras by generalizing the formalism of the bicovariant bimodule theory of 
Ref.[10] to our bimodule presented in the previous section. Since $\Pi^{n}$ and $\Theta^{n}_{~m}$ form a right-invariant basis, we are in the case where we can apply the different stages of Theorem (2,3) of [10] to state that there exist linear functionals $f^{n}_{~m}$, $f^{nk}_{~m}$, $f^{n}_{km}$ and $f^{n\ell}_{mk}$ $\in {\cal B}'$ such that\\\begin{eqnarray} \Pi^{n} a &=& (a \star f^{n}_{~k}) \Pi^{k} + (a \star f^{n\ell}_{~k}) \Theta^{k}_{~\ell} ,\\\Theta^{n}_{~m} a &=& (a \star f^{n}_{mk}) \Pi^{k} + (a \star f^{n\ell}_{mk}) \Theta^{k}_{~\ell}\end{eqnarray}and\\\begin{eqnarray*}b \Pi^{n} &=& \Pi^{k}(b \star f^{n}_{~k} \circ S) + \Theta^{k}_{~\ell}(b \star f^{n\ell}_{~k} \circ S) ,\\b \Theta^{n}_{~m} &=& \Pi^{k} (b \star f^{n}_{mk} \circ S) + \Theta^{k}_{~\ell}(b \star f^{n\ell}_{mk} \circ S)\end{eqnarray*}
where the convolution product of a functional  $f \in {\cal B}'$ and an element $a$ of ${\cal B}$ is defined as $a \star f = (f \otimes id)\Delta(a)$. From (10)and (11), one deduce 
\begin{eqnarray}
{s}(I)= \delta^{n}_{k} &,&~~f^{n\ell}_{mk}(I) = \delta^{n}_{k} \delta^{\ell}_{m},\nonumber\\f^{n}_{mk}(I) = 0 &and&~~f^{n\ell}_{~k}(I)=0
\end{eqnarray}
by setting $a=I$, and\\
\begin{eqnarray*}\Pi^{n}ab&=&(ab \star f^{n}_{~k}) \Pi^{k}+(ab \star f^{n\ell}_{~k}) \Theta^{k}_{~\ell} \nonumber\\ =((a \star f^{n}_{~\ell})(b \star f^{\ell}_{~k})&+&(a \star f^{np}_{~q})(b \star f^{q}_{pk})) \Pi^{k} +((a \star f^{n}_{~q})(b \star f^{q\ell}_{~k})\\&+&(a \star f^{np}_{~q})(b \star f^{q\ell}_{pk})) \Theta^{k}_{~\ell} ,\\\Theta^{n}_{~m}ab&=&(ab \star f^{n}_{mk}) \Pi^{k}+(ab \star f^{n\ell}_{mk}) \Theta^{k}_{\ell} \nonumber\\=((a \star f^{n}_{mq})(b \star f^{q}_{~k})&+&(a \star f^{nq}_{mp})(b \star f^{p}_{qk})) \Pi^{k}+((a \star f^{n}_{mq})(b \star f^{q\ell}_{~k}) \\&+& (a \star f^{np}_{mq})(b \star f^{q\ell}_{pk})) \Theta^{k}_{~\ell}\end{eqnarray*}for any $a,b \in {\cal B}$. Comparing the coefficients multiplying $\Pi$ and $\Theta$ and then applying $\varepsilon$, we get\\
\begin{eqnarray}
f^{n}_{~m}(ab)&=&f^{n}_{~k}(a)f^{k}_{~m}(b)+f^{nk}_{~\ell}(a)f^{\ell}_{km}(b) ,\\f^{nk}_{m\ell}(ab)&=&f^{np}_{mq}(a)f^{qk}_{p\ell}(b)+f^{n}_{mq}(a)f^{qk}_{~l}(b) ,\\f^{n}_{m\ell}(ab)&=&f^{np}_{mq}(a)f^{q}_{p\ell}(b)+f^{n}_{mk}(a)f^{k}_{~\ell}(b) ,\\f^{nk}_{~\ell}(ab)&=&f^{n}_{~q}(a)f^{qk}_{~\ell}(b)+f^{np}_{~q}(a)f^{qk}_{p\ell}(b).\end{eqnarray}Applying now $\Delta_{L}$ on both sides of (10) and (11), we obtain respectively\\\begin{eqnarray}\Lambda^{n}_{~p} a_{(1)} \otimes \Pi^{p} a_{(2)}&+&\Lambda^{n}_{~p} S(p^{q})a_{(1)} \otimes \Theta^{p}_{~q} a_{(2)}=f^{n}_{~p}(a_{(1)})a_{(2)} \Lambda^{p}_{~k} \otimes a_{(3)} \Pi^{k} \nonumber \\ +(f^{nq}_{~p}(a_{(1)}) a_{(2)} \Lambda^{p}_{~k} S( \Lambda^{\ell}_{~q})&+&f^{n}_{~p}(a_{(1)}) a_{(2)}  \Lambda^{p}_{~k} S(p^{\ell})) \otimes a_{(3)} \Theta^{k}_{~\ell}\end{eqnarray}and\begin{eqnarray}\Lambda^{n}_{~q} S( \Lambda^{p}_{~m}) a_{(1)} \otimes \Theta^{q}_{~p} a_{(2)} &=&f^{n}_{mp}(a_{(1)})a_{(2)} \Lambda^{p}_{~k} \otimes a_{(3)} \Pi^{k} \nonumber\\
 \Lambda^{p}_{~k} S(p^{\ell})&+&f^{nq}_{mp}(a_{(1)}) a_{(2)} \Lambda^{p}_{~k} S( \Lambda^{\ell}_{~q})) \otimes a_{(3)} \Theta^{k}_{~\ell}.
\end{eqnarray}
By virtue of (10) and (11), the first hand sides of previous equations can berewritten respectively as:\\\begin{eqnarray*}(\Lambda^{n}_{~p} a_{(1)}f^{p}_{~k}(a_{(2)})&+&\Lambda^{n}_{~p} S(p^{q})f^{p}_{qk}(a_{(2)})) \otimes a_{(3)} \Pi^{k}+( \Lambda^{n}_{~p} a_{(1)}f^{p\ell}_{~k}(a_{(2)})\\ &+& \Lambda^{n}_{~p} S(p^{q})a_{(1)}f^{p\ell}_{qk}(a_{(2)})) \otimes a_{(3)}\Theta^{k}_{~\ell}\end{eqnarray*}and \\\begin{eqnarray*}\Lambda^{n}_{~q}S( \Lambda^{p}_{~m})a_{(1)}f^{q}_{pk}(a_{(2)}) \otimes a_{(3)} \Pi^{k}+\Lambda^{n}_{~q} S( \Lambda^{p}_{~m}) a_{(1)}f^{q\ell}_{pk}(a_{(2)}) \otimes a_{(3)} \Theta^{k}_{~\ell}.\end{eqnarray*}Substituting the first expression into the first hand side of (17) and the second into the first hand side of (18), comparing the coefficients multiplying $I \otimes \Pi$ and $I \otimes \Theta$
 and then applying $ I \otimes \varepsilon$, we obtain\\\begin{eqnarray}\Lambda^{n}_{~k}(f^{k}_{~m} \star a)+\Lambda^{n}_{~k}S(p^{\ell})(f^{k}_{\ell m} \star a)&=&(a \star f^{n}_{~k}) \Lambda^{k}_{~m}\\\Lambda^{n}_{~k} (f^{k\ell}_{~m} \star a)+\Lambda^{n}_{~k} S(p^{q})(f^{k\ell}_{qm} \star a)&=&(a \star f^{n}_{~k}) \Lambda^{k}_{~m} S(p^{\ell})+(a \star f^{nk}_{~q}) \Lambda^{q}_{~m} S( \Lambda^{\ell}_{~k})\nonumber\\\\\Lambda^{n}_{~k} S( \Lambda^{q}_{~m})(f^{k}_{q\ell} \star a)&=&(a \star f^{n}_{mk}) \Lambda^{k}_{~\ell}\\\Lambda^{n}_{~k} S( \Lambda^{q}_{~m})(f^{kp}_{ql} \star a)&=&(a \star f^{n}_{mk}) \Lambda^{k}_{~\ell} S(p^{p})+(a \star f^{nq}_{mk}) \Lambda^{k}_{~\ell} S( \Lambda^{p}_{~q})\nonumber\\
\end{eqnarray}
for any $a \in {\cal B}$. In the following, we assume that\\
\begin{eqnarray}
f^{nk}_{m\ell} = \tilde{f}^{k}_{~m} \star f^{n}_{\ell}  &~~f^{nk}_{~m} = \tilde{\eta}^{k} \star f^{n}_{~m}  &and~~f^{n}_{km} = \tilde{\eta}_{k} \star f^{n}_{~m}\end{eqnarray}
where $\tilde{\eta}^{n}$, $f^{n}_{~m} \in {\cal B}'$  and the convolution product of two functionals $ \in {\cal B}'$ is defined as $(f_{1} \star f_{2})(a) = (f_{1} \otimes f_{2}) \Delta(a)$ for any $a \in {\cal B}$. The substitution of (23) into (12) gives\\
\begin{eqnarray}
\tilde{\eta}^{n}(I) = 0 &,~~ \tilde{\eta}_{n}(I) = 0  &and~~ \tilde{f}^{n}_{~m}(I) = \delta^{n}_{m}
\end{eqnarray}
and the substitution of ( 23) into (13-16) and (19) gives respectively\\
\begin{eqnarray}
f^{n}_{~m}(ab)=(\varepsilon(a_{(1)})\varepsilon(b_{(1)}) + \tilde{\eta}^{k} (a_{(1)})\tilde{\eta}_{k}(b_{(1)}))f^{n}_{~\ell}(a_{(2)})f^{\ell}_{~m}(b_{(2)}),
\end{eqnarray}
\begin{eqnarray}
f^{nk}_{m\ell}(ab)&=&\tilde{f}^{k}_{~m}(a_{(1)}b_{(1)})f^{n}_{~\ell}(a_{(2)}b_{(2)}) = \nonumber \\ 
&=&(\tilde{f}^{p}_{~m}(a_{(1)})\tilde{f}^{k}_{~p}(b_{(1)})+\tilde{\eta}_{m}(a_{(1)})\tilde{\eta}^{k}(b_{(1)}))f^{n}_{~q}(a_{(2)})f^{q}_{~\ell}(b_{(2)}) ,
\end{eqnarray}
\begin{eqnarray}
f^{n}_{km}(ab)&=&\tilde{\eta}_{k}(a_{(1)}b_{(1)})f^{n}_{~m}(a_{(2)}b_{(2)})=\nonumber\\
&=&(\tilde{f}^{p}_{~k}(a_{(1)}) \tilde{\eta}_{p}(b_{(1)})+\tilde{\eta}_{k}(a_{(1)}) \varepsilon (b_{(1)}))f^{n}_{~q}(a_{(2)})f^{q}_{~m}(b_{(2)}),
\end{eqnarray}
\begin{eqnarray}
f^{nk}_{~m}(ab)&=&\tilde{\eta}^{k}(a_{(1)}b_{(1)})f^{n}_{~m}(a_{(2)}b_{(2)})=\nonumber \\
&=&(\varepsilon(a_{(1)})\tilde{\eta}^{k}(b_{(1)})+\tilde{\eta}^{p}(a_{(1)})\tilde{f}^{k}_{p}(b_{(1)}))f^{n}_{~q}(a_{(2)})f^{q}_{~m}(b_{(2)})
\end{eqnarray}
and \\ 
\begin{eqnarray}
\Lambda^{n}_{~k}(f^{k}_{~m} \star a) + \Lambda^{n}_{~k}S(p^{l})(\tilde{\eta}_{l} \star f^{k}_{~m} \star a) = (a \star f^{n}_{~k}) \Lambda^{k}_{~m}.
\end{eqnarray}
Using now the associativity of the convolution product and making the substitution of (29) into the second hand of (20), (21) and (22), we get\\ 
\begin{eqnarray}
\tilde{\eta}^{i} \star f^{k}_{~j} \star a &+& S(p^{\ell})(\tilde{f}^{i}_{~\ell}\star f^{k}_{~j} \star a)=(f^{k}_{~j} \star a)S(p^{i}) + (f^{k}_{~j} \star a \star \tilde{\eta}^{l})S( \Lambda^{i}_{~l})\nonumber\\
&+& S(p^{\ell})(\tilde{\eta}_{\ell} \star f^{k}_{~j} \star a)S(p^{i})+S(p^{p})(\tilde{\eta}_{p} \star \tilde{\eta}^{\ell})S(\Lambda^{i}_{~\ell}),\nonumber\\
S(\Lambda^{k}_{~m})( \tilde{\eta}_{k} \star f^{n}_{~\ell} \star a)&=&(f^{n}_{~\ell} \star a \star \tilde{\eta}_{m})+S(p^{k})(\tilde{\eta}_{k} \star f^{n}_{\ell} \star a \star \tilde{\eta}_{m})
\end{eqnarray}
and\\
\begin{eqnarray}
S( \Lambda^{\ell}_{~m})( \tilde{f}^{i}_{~\ell} \star f^{k}_{~j} \star a)&=&S( \Lambda^{\ell}_{~m})( \tilde{\eta}_{\ell} \star f^{k}_{~j} \star a)S(p^{i})+(f^{k}_{~j} \star a \star \tilde{f}^{\ell}_{~m})S( \Lambda^{i}_{~\ell})\nonumber\\&+&S(p^{q})(\tilde{\eta}_{q} \star f^{k}_{~j} \star a \star \tilde{f}^{\ell}_{~m})S( \Lambda^{i}_{~\ell})\end{eqnarray}where we have multiplied from the left by $S( \Lambda^{p}_{n}))$ and we have used (21) before (29) to obtain (32). We see from these relations that for any $a \in {\cal A}$, (30) is a commutation rule between elements of ${\cal A}$ and $p^{n}$ and (29), (31) and (32) involve $p^{n}$ in the commutation rules between elements of ${\cal A}$ . Demanding ${\cal A}$ to be a Hopf subalgebra of ${\cal B}$ (condition 2.) leads to\\{\bf Proposition (3,1)}: ${\cal A}$ is a Hopf sub-algebra of ${\cal B}$ iff \\\begin{eqnarray}\tilde{\eta}_{n}(a) = 0 &,~~ a \in {\cal A} \end{eqnarray}$proof.$ For ${\cal A}$ to be a Hopf subalgebra of ${\cal B}$, (29) and (32)must involve commutation rules between elements of ${\cal A}$ only. this is satisfied if the term $\Lambda^{n}_{~k}S(p^{\ell})(\tilde{\eta}_{\ell} \star f^{k}_{~m} \star a) $ of (29) vanishes. Multiplying it from the left by $S(\Lambda^{i}_{~n})$ and using the fact that $S(p^{\ell})$ are linearly independent and generate with the unity the bimodule $\Gamma_{p}$ (see condition 5.), one obtains $\tilde{\eta}_{\ell} \star f^{i}_{~m} \star a = 0$ which permits us to write, for any $a$ and $b \in {\cal A}$, (25) under the form\\
\begin{eqnarray}
f^{n}_{~m}(ab) = f^{n}_{~k}(a)f^{k}_{~m}(b) \Rightarrow f^{n}_{~k} \star f^{k}_{~m} \circ S=f^{n}_{~k} \circ S \star f^{k}_{~m}=\delta^{n}_{m}
\end{eqnarray}
Replacing $ a $ by $f^{m}_{~k} \circ S \star a$ into $\tilde{\eta}_{\ell} \star f^{i}_{~m} \star a=0$ and then acting $\varepsilon$, we get $\tilde{\eta}_{\ell}(a) = 0$ for any $a \in {\cal A}$ which implies also that equation (31) is trivial and the term of (32) which involves $p^{n}$ in the commutation rules between elements of ${\cal A}$ vanishes. Conversely, it is easy to see that if $\tilde{\eta}_{\ell}(a) = 0$ for any $a \in {\cal A}$, the relations (29) and (32) involve commutation rules between elements of ${\cal A}$ only and, therefore, ${\cal A}$ is a Hopf subalgebra of ${\cal B}$. Q.ED.\\By virtue of this proposition, one deduces certain results which will be used in the following.\\1. (11) can be written as: \begin{eqnarray*}\Theta^{n}_{~m}a = (a \star f^{n\ell}_{mk})\Theta^{k}_{~\ell},  &a \in {\cal A} \end{eqnarray*}which shows that the set of the right invariant one-forms $\Theta^{n}_{~m}$ generates a subimodule, $\Gamma_{\cal A} \subset \Gamma$, over 
${\cal A}$ and for $\tilde{\eta}_{n}(p^{k})=0$ it generates a ${\cal B}$-subimodule of $\Gamma$. since it transforms according to an adjoint representation $\Lambda$ of ${\cal H}$ (6), the assumption $f^{n\ell}_{mk} = \tilde{f}^{\ell}_{~m} \star f^{n}_{~k}$ is justified [11].\\2. For any $a \in {\cal A}$, (26) can be written as:\\\begin{eqnarray*}\tilde{f}^{k}_{~m}(a_{(1)}b_{(1)})f^{n}_{~\ell}(a_{(2)}b_{(2)})&=&\tilde{f}^{p}_{~m}(a_{(1)}) \tilde{f}^{k}_{~p}(b_{(1)})f^{n}_{~q}(a_{(2)})f^{q}_{~\ell}(b_{(2)})\nonumber\\&=&\tilde{f}^{p}_{~m}(a_{(1)}) \tilde{f}^{k}_{~p}(b_{(1)})f^{n}_{~\ell}(a_{(2)}b_{(2)})\end{eqnarray*}where we have used (34) to get the third hand side. Multiplying both sides from the right by $f^{\ell}_{~r}(S(a_{(3)}b_{(3)}))$ and using again (34), we get\\
\begin{eqnarray}
\tilde{f}^{n}_{~m}(ab)= \tilde{f}^{k}_{~m}(a) \tilde{f}^{n}_{~k}(b),    &a,b \in {\cal A} .
\end{eqnarray}
Following similar considerations, one obtains from (28)\\
\begin{eqnarray}
\tilde{\eta}^{n}(ab)= \varepsilon + \tilde{\eta}^{m}(a) \tilde{f}^{n}_{~m}(b),  &a,b \in {\cal A} 
\end{eqnarray}
which leads to\\
\begin{eqnarray}
\tilde{\eta}^{\ell}( \Lambda^{n}_{~k}S( \Lambda^{k}_{~m})) &=& 
\tilde{\eta}^{\ell}( \delta^{n}_{m}) = 0 = \tilde{\eta}^{\ell}
(S(\Lambda^{n}_{~m})) + \tilde{\eta}^{q}( \Lambda^{n}_{~k})
\tilde{f}^{\ell}_{~q}(S( \Lambda^{k}_{~m})),\nonumber\\\tilde{\eta}^{\ell}
(S( \Lambda^{n}_{~k}) \Lambda^{k}_{~m})&=&\tilde{\eta}^{\ell}(\delta^{n}_{m}) = 0 = \tilde{\eta}^{\ell}( \Lambda^{n}_{~m}) + \tilde{\eta}^{q}
(S( \Lambda^{n}_{~k}))\tilde{f}^{\ell}_{q}( \Lambda^{k}_{~m}). 
\end{eqnarray}
Now replacing $ab$ by $-S(\Lambda^{a}_{~b})p^{b}=S(p^{a})$ into (28) and 
using (37), one obtains 
\begin{eqnarray}
\tilde{\eta}^{n}(S(p^{a}))=- \tilde{\eta}^{n}(p^{a})- \tilde{\eta}^{k}
(S(\Lambda^{a}_{~b}))\tilde{f}^{n}_{~k}(p^{b}).
\end{eqnarray}
In the other hand, by using (33), (27) gives \\
\begin{eqnarray*}
(\tilde{\eta}_{k} \star f^{n}_{~m})(S(p^{a}))&=&\tilde{\eta}_{k}(S(p^{c}))
f^{n}_{~m}(S( \Lambda^{a}_{~c})) = -(\tilde{\eta}_{k} \star f^{n}_{m})
(S( \Lambda^{a}_{~b})p^{b})\nonumber\\
&=& - \tilde{f}^{p}_{~k}(S( \Lambda^{c}_{~b}))\tilde{\eta}_{p}(p^{b})
f^{n}_{~m}(S( \Lambda^{a}_{~c}))
\end{eqnarray*}
which can be multiplied from the right by $f^{m}_{~\ell}( \Lambda^{c}_{~d})$ to give\\
\begin{eqnarray}
\tilde{\eta}_{n}(S(p^{a})) = - \tilde{f}^{k}_{~n}(S( \Lambda^{a}_{~b})) 
\tilde{\eta}_{k}(p^{b}).
\end{eqnarray}
A similar computation gives respectively from (25) and (26)\\ 
\begin{eqnarray}
f^{n}_{~m}(p^{a})=-f^{n}_{~m}( \Lambda^{a}_{~b}S(p^{b})) = -f^{n}_{~k}
(\Lambda^{a}_{~b})f^{k}_{~m}(S(p^{b})) - \tilde{\eta}^{k}( \Lambda^{a}_{~b}) 
\tilde{\eta}_{k}(S(p^{b})) \delta^{n}_{m}
\end{eqnarray}
and
\begin{eqnarray}
\tilde{f}^{n}_{~m}(p^{a})=-\tilde{f}^{k}_{~m}(\Lambda^{a}_{~b})
\tilde{f}^{n}_{~k}(S(p^{b})) +\tilde{\eta}^{k}(\Lambda^{c}_{~b})\tilde{\eta}_{k}(S(p^{b}))\tilde{f}^{n}_{~m}(\Lambda^{a}_{~c}).
\end{eqnarray}
We are now ready to investigate the different commutation rules 
bra ${\cal B}$.\\ First, due to the proposition (3,1), (29) reduces to\\
\begin{eqnarray}
\Lambda^{n}_{~k}(f^{k}_{~m} \star a) =(a \star f^{n}_{~k}) \Lambda^{k}_{~m},
\end{eqnarray}
for any $a \in {\cal A}$ and gives
\begin{eqnarray}
\Lambda^{n}_{~m}p^{a} = f^{n}_{~k}( \Lambda^{a}_{~b})p^{b} \Lambda^{k}_{~m} + 
f^{n}_{~k}(p^{a}) \Lambda^{k}_{~m} - \Lambda^{n}_{~k} \Lambda^{a}_{~b} 
f^{k}_{~m}(p^{b}) - \Lambda^{n}_{~m}S(p^{\ell}) \Lambda^{a}_{~b} 
\tilde{\eta}_{\ell}(p^{b})\end{eqnarray}
for $a=p^{a}$.\\
Secondly, replacing $a \in {\cal B}$ by $S(a)$ into (30 - 32) and then acting 
on both sides $S^{-1}$, we get respectively\\
\begin{eqnarray}
(a \star f^{k}_{~m} \circ S \star \tilde{\eta}^{\ell} \circ S)&+&(a \star 
f^{k}_{~m} \circ S \star \tilde{f}^{\ell}_{~q} \circ S)p^{q} = p^{\ell}
(a \star f^{k}_{~m} \circ S)\nonumber\\
&+&p^{\ell}(a \star f^{k}_{~m} \circ S \star \tilde{\eta}_{q} \circ S)p^{q}+
\Lambda^{\ell}_{~p}(\tilde{\eta}^{p}\circ S \star a \star f^{k}_{~m} \circ S)
\nonumber\\
&+&\Lambda^{\ell}_{~p}(\tilde{\eta}^{p} \circ S \star a \star f^{k}_{~m} 
\circ S \star \tilde{\eta}_{q} \circ S)p^{q},\\(a \star f^{n}_{~\ell} \circ S 
\star \tilde{\eta}_{k} \circ S)\Lambda^{k}_{~m}
&=&\tilde{\eta}_{m} \circ S \star a \star f^{n}_{~\ell} \circ S+
(\tilde{\eta}_{m} \circ S \star a \star f^{n}_{~\ell} \circ S \star 
\tilde{\eta}_{q} \circ S)p^{q}\nonumber\\
\end{eqnarray}
and\\
\begin{eqnarray}(a \star f^{k}_{~j} \circ S \star \tilde{f}^{i}_{\ell} \circ S) \Lambda^{\ell}_{~m}&=&p^{i}(a \star f^{k}_{~j} \circ S \star \tilde{\eta}_{\ell} \circ S) \Lambda ^{\ell}_{~m}+\Lambda^{i}_{\ell}(\tilde{f}^{\ell}_{~m} \circ S \star a \star f^{k}_{~j} \circ S)\nonumber\\&+&\Lambda^{i}_{\ell}
(\tilde{f}^{\ell}_{~m} \circ S \star a \star f^{k}_{~j} \circ S \star 
\tilde{\eta}_{q} \circ S)p^{q}.
\end{eqnarray}
Due to the proposition (3,1), (45) is trivial for $a \in {\cal A}$ and by 
replacing $a$ by $(a \star f^{n}_{~k}) \in {\cal A}$, the equations (46) and 
(44) give respectively\\
\begin{eqnarray}
e{f}^{k}_{~m} \circ S \star a)=(a \star \tilde{f}^{n}_{~k} \circ S) 
\Lambda^{k}_{~m} 
\end{eqnarray}
and\\
\begin{eqnarray}
p^{n}a=(a \star \tilde{f}^{n}_{~k} \circ S)p^{k} + a \star \tilde{\eta}^{n} 
\circ S - \Lambda^{n}_{~k}(\tilde{\eta}^{k} \circ S \star a).
\end{eqnarray}
For $a = p^{a}$, the equations (45) and (46) give respectively\\
\begin{eqnarray}
\Lambda^{k}_{~m} \tilde{\eta}_{k}(S(p^{a}))=\tilde{\eta}_{m}(S(p^{b})) 
\Lambda^{a}_{~b}
\end{eqnarray}
and \\
\begin{eqnarray}
\Lambda^{n}_{~m}p^{b}( \delta^{a}_{b} + \tilde{\eta}_{b}(S(p^{a}))) &=& 
(\tilde{f}^{n}_{~k}(S( \Lambda^{a}_{~b})) - \delta^{n}_{~b} \tilde{\eta}_{k}
(S(p^{a})))p^{b} \Lambda^{k}_{~m}\nonumber\\ 
&+& \tilde{f}^{n}_{~k}(S(p^{a})) \Lambda^{k}_{~m} - \Lambda^{n}_{~k} \Lambda^{a}_{b} \tilde{f}^{k}_{~m}(S(p^{b})).
\end{eqnarray}
To compare (50) with (43), we have to compute the fourth term of the second 
hand of (43) which can be written as\\
\begin{eqnarray}
- \Lambda^{n}_{~m}S(p^{\ell}) \Lambda^{a}_{~c} \tilde{\eta}_{\ell}(p^{c}) &=& 
\Lambda^{n}_{~m}S( \Lambda^{\ell}_{~k})p^{k} \Lambda^{a}_{~c} 
\tilde{\eta}_{\ell}(p^{c})\nonumber\\
= - \Lambda^{n}_{~m}p^{\ell}\tilde{\eta}_{\ell}(p^{a}) &+& \Lambda^{n}_{~m} 
\tilde{\eta}^{k}(S( \Lambda^{a}_{~b})) \tilde{\eta}_{k}(p^{b}) - 
\Lambda^{n}_{~m} \Lambda^{a}_{~b} \tilde{\eta}^{k}(S( \Lambda^{b}_{~c})) 
\tilde{\eta}_{k}(p^{c})\nonumber\\
\end{eqnarray}
where we have used (39),\\
\begin{eqnarray}
p^{k} \Lambda^{a}_{~c} = \tilde{f}^{k}_{~\ell}(S( \Lambda^{a}_{~b})) 
\Lambda^{b}_{~c}p^{\ell} + \tilde{\eta}^{k}(S( \Lambda^{a}_{~b})) 
\Lambda^{b}_{~c} - \Lambda^{k}_{~\ell} \Lambda^{a}_{~b} \tilde{\eta}^{\ell}
(S( \Lambda^{b}_{~c}))
\end{eqnarray}
obtained from (48) by setting $a= \Lambda^{a}_{~c}$ and\\
\begin{eqnarray}
S( \Lambda^{k}_{~m}) \Lambda^{a}_{~c} \tilde{\eta}_{k}(p^{c})= 
\tilde{\eta}_{m}(p^{a})
\end{eqnarray}
obtained from (21) by setting $a=p^{a}$ and by using (23) and the proposition 
(3,1).\\
Replacing the fourth term of the right hand side of (43) by (51), we get:\\
\begin{eqnarray}
\Lambda^{n}_{~m}p^{b}(\delta^{a}_{b} + \tilde{\eta}_{b}(S(p^{a}))) &=& 
f^{n}_{~k}( \Lambda^{a}_{~b})p^{b} \Lambda^{k}_{~m} + (f^{n}_{~k}(p^{a}) + 
\delta^{n}_{k} \tilde{\eta}^{p}(S( \Lambda^{a}_{~b}))\tilde{\eta}_{p}(p^{b})) 
\Lambda^{k}_{~m}\nonumber\\ &-& \Lambda^{n}_{~k} \Lambda^{a}_{~b}
(f^{k}_{~m}(p^{b}) + \delta^{k}_{m} \tilde{\eta}^{p}(S( \Lambda^{b}_{~c}))
\tilde{\eta}_{p}(p^{c}))
\end{eqnarray}
which can be written, after the make of use of (37), (39) and (40), under the 
form\\
\begin{eqnarray}
\Lambda^{n}_{~m}p^{b}(\delta^{a}_{b} &+& \tilde{\eta}_{b}(S(p^{a}))) = 
f^{n}_{~k}(\Lambda^{a}_{~c})p^{c}\Lambda^{k}_{~m}\nonumber\\ &-& f^{n}_{~\ell}
( \Lambda^{a}_{~b})f^{\ell}_{~k}(S(p^{b})) \Lambda^{k}_{~m} + \Lambda^{n}_{~k}  \Lambda^{a}_{~b}f^{k}_{~\ell}( \Lambda^{b}_{~c})f^{\ell}_{~m}(S(p^{c})).
\end{eqnarray}
Therefore, by virtue of the condition 5. of the inhomogeneous Hopf algebra 
definition, we can deduce from (55) and (50) that\\
\begin{eqnarray}
\tilde{f}^{n}_{~k}(S( \Lambda^{a}
_{~b})) = f^{n}_{~k}( \Lambda^{a}_{~b}) + \delta^{n}_{b} \tilde{\eta}_{k}
(S(p^{a}))\end{eqnarray}and\\\begin{eqnarray}\tilde{f}^{n}_{~k}(S(p^{a})) = 
-f^{n}_{~\ell}( \Lambda^{a}_{b})f^{\ell}_{~k}(S(p^{b})) = f^{n}_{~k}(p^{a}) + 
\delta^{n}_{k} \tilde{\eta}^{\ell}(S( \Lambda^{a}_{~b})) \tilde{\eta}_{\ell}
(p^{b}).
\end{eqnarray}
where we have used (40) to get the second relation. Let us note that, for 
$a= \Lambda^{a}_{b}$, the equations (47) and (42) are consistent in virtue of 
(56) and (49).\\ Now, multiplying both sides of (55) from the left by 
$f(S(\Lambda))$, we get\begin{eqnarray}p^{n} \Lambda^{k}_{~m} &=& f^{k}_{~q}
(S( \Lambda^{n}_{~p}))( \delta^{p}_{\ell} + \tilde{\eta}_{\ell}(S(p^{p}))) 
\Lambda^{q}_{~m}p^{\ell}\nonumber\\ &+& f^{k}_{~q}(S(p^{n})) \Lambda^{q}_{~m} - \Lambda^{n}_{~q} \Lambda^{k}_{~p}f^{p}_{~m}(S(p^{q}))
\end{eqnarray}
where we have used (42). Comparing the latter equation with (52), we obtain\\
\begin{eqnarray}
\tilde{f}^{n}_{~\ell}(S( \Lambda^{k}_{~q})) = f^{k}_{~q}(S( \Lambda_{l} + 
\tilde{\eta}_{\ell}(S(p^{p})))
\end{eqnarray}
and\\
\begin{eqnarray}\tilde{\eta}^{n}(S( \Lambda^{k}_{~q})) = f^{k}_{~q}(S(p^{n})).
\end{eqnarray}
Substituting (56) into (59), we get\\
\begin{eqnarray}
f^{n}_{~\ell}(\Lambda^{k}_{~q}) + \delta^{n}_{q} \tilde{\eta}_{\ell}(S(p^{k})) = f^{k}_{~q}(S( \Lambda^{n}_{~p}))(\delta^{p}_{\ell} + \tilde{\eta}^{\ell}
(S(p^{p})))
\end{eqnarray}
which can be written, after multiplying both sides by $f^{m}_{k}
(\Lambda^{t}_{~n})$ as\\\begin{eqnarray}(f^{m}_{~k}( \Lambda^{t}_{~n}) - 
\delta^{m}_{~n} \delta^{t}_{~k})(f^{n}_{~\ell}(\Lambda^{k}_{~q}) + 
\delta^{n}_{q} \delta^{k}_{\ell} + \delta^{n}_{~q} \tilde{\eta}_{\ell}(S(p^{k})))=0
\end{eqnarray}                 
or, using again (56),\\
\begin{eqnarray}
(\tilde{f}^{m}_{~k}(S( \Lambda^{t}_{~n})) - \delta^{m}_{n}\delta^{t}_{k} - 
\delta^{m}_{n}\tilde{\eta}_{k}(S(p^{t})))(\tilde{f}^{n}_{~\ell}
(S( \Lambda^{k}_{~q}))+\delta^{n}_{q}\delta^{k}_{\ell})=0.
\end{eqnarray}
which are the characteristic
 equations for the matrices $R^{mt}_{nk}= f^{m}_{~k}(\Lambda^{t}_{~n})$ and 
$\tilde{R}^{mt}_{nk}= \tilde{f}^{m}_{~k}(S(\Lambda^{t}_{~n}))$. To have the 
commutation rules between the translations, it suffices to replace $a$ by 
$p^{a}$ into (44) to get after some straightforward computation\\
\begin{eqnarray}
p^{n}p^{m} &=& (\tilde{f}^{n}_{~k}(S( \Lambda^{m}_{~\ell})) - \delta^{n}_{\ell} \tilde{\eta}_{k}(S(p^{m})))p^{\ell}p^{k} + \tilde{\eta}^{n}(S( \Lambda^{m}_{~k}))p^{k} + \tilde{f}^{n}_{~k}(S(p^{m}))p^{k}\nonumber\\
&+&\tilde{\eta}^{n}(S(p^{m})) - \Lambda^{n}_{~k} \Lambda^{m}_{\ell} 
\tilde{\eta}^{k}(S(p^{\ell}))\nonumber\\&=& f^{n}_{~k}(\Lambda^{m}_{~\ell})
p^{\ell}p^{k} - (f^{n}_{~k}(\Lambda^{m}_{~\ell})-\delta^{n}_{~\ell}
\delta^{m}_{~k})\tilde{\eta}^{\ell}(S( \Lambda^{k}_{~p}))p^{p} \nonumber\\
&+&\tilde{\eta}^{n}(S(p^{m})) - \Lambda^{n}_{~k} \Lambda^{m}_{\ell} 
\tilde{\eta}^{k}(S(p^{\ell}))
\end{eqnarray}
where we have used (56), (57) and (60). With the assumption (23), it follows from the previous theorem (3,1):Let ${\cal B}$ an inhomogeneous Hopf algebra satisfying 1.-5. on which there exists a bicovariant differential calculus, Then there exist functionals $f^{n}_{~m}$, $\tilde{f}^{n}_{~m}$, $\tilde{\eta}^{n}$ and 
$\tilde{\eta}_{m} \in {\cal B}'$ satisfying (24), (25-28), (33), (56-57) and 
(59-60) such that the relations (42), (48) and (64) are satisfied.\\
\section{\bf Consistency Conditions}Here we shall continue the study of the 
consistency conditions between the different functionals and the commutation 
rules of the generators of ${\cal G}$. In the following we set 
$R^{kn}_{m\ell}=f^{k}_{~\ell}(\Lambda^{n}_{~m})$, $\tilde{R}^{kn}_{m\ell}=
\tilde{f}^{k}_{~\ell}(S(\Lambda^{n}_{~m}))$, $Q^{n}_{k}=\tilde{\eta}_{k}
(S(p^{n}))$, $Z^{nm}_{k}=\tilde{\eta}^{n}(S(\Lambda^{m}_{~k})) = f^{m}_{~k}(S(p^{n}))$  (see (60)) , $\tilde{Z}^{nm}_{k}=\tilde{f}^{n}_{k}(S(p^{m})) = 
-R^{nm}_{pq}Z^{pq}_{k}$(see (57)) and $T^{nm}=\tilde{\eta}^{n}(S(p^{m}))$.
 We start this study by  \\{\bf L
emma (4,1)}\\
\begin{eqnarray}
Q^{n}_{m}=\lambda\delta^{n}_{m}    &,~~\lambda \in {\cal C},\lambda \not= -1
\end{eqnarray}
$proof.$ Applying $f^{a}_{~b}$ and $\tilde{f}^{a}_{~b} \circ S$ on both sides of (49), we obtain respectively $(I \otimes Q)R = R(Q \otimes I)$ and 
$(I \otimes Q)\tilde{R} = \tilde{R}(Q \otimes I)$ implying, due to (56) 
($\tilde{R}=R + I \otimes Q$), $I \otimes Q^{2} = Q \otimes Q$ which solution is of the form (65). For $\lambda = -1$, we see from (62) that the $R$-matrix is 
not invertible.    Q.E.D.\\Applying $f^{a}_{~b}$ on both sides of (42) and 
$\tilde{f}^{a}_{~b} \circ S$on both sides of (47) for $a=\Lambda^{p}_{~q}$ and 
using (34-35), we get the Yang-Baxter equations for the matrices $R$ 
\begin{eqnarray}
(I \otimes R)(R\otimes I)(I \otimes R)=(R \otimes I)(I \otimes R)(R \otimes I)
\end{eqnarray}
and $\tilde{R}$
\begin{eqnarray}
(I \otimes \tilde{R})(\tilde{R} \otimes I)(I \otimes \tilde{R})=(\tilde{R} 
\otimes I)(I \otimes \tilde{R})(\tilde{R} \otimes I).
\end{eqnarray}
s tied to the $R$-matrix by (56), it is necessary to check the consistency 
between both Yang-Baxter equations. Applying $f^{a}_{~b}$ on both sides of (47) and $\tilde{f}^{a}_{~b} \circ S$ on both sides of (42) for 
$a = \Lambda^{p}_{~q}$, we obtain \begin{eqnarray}(R \otimes I)(I \otimes R)
(\tilde{R} \otimes I) = (I \otimes \tilde{R})(R \otimes I)(I \otimes R)
\end{eqnarray}
and
\begin{eqnarray}
(\tilde{R} \otimes I)(I \otimes \tilde{R})(R \otimes I ) = (I \otimes R)
(\tilde{R} \otimes I)(I \otimes \tilde{R}).\end{eqnarray} Replacing (56) into 
both sides of (68)  and (69) and using (66), we obtain respectively 
\begin{eqnarray*}
(R \otimes I)(I \otimes R)(I \otimes Q \otimes I) &=& (I \otimes I \otimes Q)
(R \otimes I)(I \otimes R)\end{eqnarray*}and \begin{eqnarray*}(R \otimes I)
(I \otimes I \otimes Q)(R \otimes I)&+&(I \otimes Q \otimes I)(I \otimes R)
(R \otimes I)\nonumber\\&+&(I \otimes Q \otimes I)(I \otimes I \otimes Q)
(R \otimes I)\nonumber\\=(I \otimes R)(R 
\otimes I)(I \otimes I \otimes Q)&+&(I \otimes R)(I \otimes Q \otimes I)
(I \otimes R)\nonumber\\&+&(I \otimes R)(I \otimes Q \otimes I)
(I \otimes I \otimes Q).\end{eqnarray*}By virtue of (65) the first equation is 
satisfied and using (62) ($R^{2}=-R(I \otimes Q) + I \otimes I + I \otimes Q$), we can see that the second equation is also satisfied . Now replacing the third factor of the left hand side and the first factor of the right hand side of (67) by (56) and using (69), we get  \begin{eqnarray*}(\tilde{R} \otimes I)
(I \otimes \tilde{R})(I \otimes Q \otimes I)=(I \otimes I \otimes Q)
(\tilde{R} \otimes I)(I \otimes \tilde{R}) \end{eqnarray*}which is satisfied in virtue of (65). Therefore (66), (67), (68) and (69) are equivalent. The action 
of $\tilde{\eta}^{a} \circ S$ on both sides of (49) gives\begin{eqnarray}ZQ = 
(I \otimes Q)Z\end{eqnarray}and on (42) and (47) gives, for 
$a= \Lambda^{p}_{~q}$,
\begin{eqnarray}
(Z \otimes I)R + (\tilde{R} \otimes I)(I \otimes Z)R = (I \otimes R)
(Z \otimes R \otimes I)(I \otimes Z)
\end
{eqnarray}
and
\begin{eqnarray}
(Z \otimes I)\tilde{R} + (\tilde{R} \otimes I)(I \otimes Z)\tilde{R} = 
(I \otimes \tilde{R})(Z \otimes I)+(I \otimes \tilde{R})(\tilde{R} \otimes I)
(I \otimes Z)
\end{eqnarray}
where we have used (36). Due to (65), 
the equation (70) is satisfied and if we use (56) into (72) and compare with 
(71), we obtain 
\begin{eqnarray*}
(Z \otimes I)(I \otimes Q) &+& 
(\tilde{R} \otimes I)(I \otimes Z)(I \otimes Q)\nonumber\\=(I \otimes I \otimes Q)(Z \otimes I) &+& (I \otimes I \otimes Q)(\tilde{R} \otimes I)(I \otimes Z)
\end{eqnarray*}
which is satisfied by virtue of (65) or (70). Therefore the relation (71) 
implies (72).\\We pass now to the action of the different functionals on the 
commutation rules between the generators of ${\cal A}$ and the translations. 
Applying $\tilde{\eta}_{a} \circ S$ on both sides of (52) and using the relations
\begin{eqnarray*}
\tilde{\eta}_{a}(aS(p^{n}))=\tilde{f}^{b}_{~a}(a)\tilde{\eta}
_{b}(S(p^{n}))\end{eqnarray*}and\begin{eqnarray*}\tilde{\eta}_{a}(S(p^{n})a)=
\tilde{\eta}_{a}(S(p^{n}))\varepsilon(a)
\end{eqnarray*}
obtained from (27) for any $a \in {\cal A}$, we get $(Q \otimes I)\tilde{R}=
\tilde{R}(I \otimes Q)$ which is satisfied due to (65). Applying $f^{a}_{~b} 
\circ S$ and $\tilde{f}^{a}_{~b} \circ S$ on both sides of (52), we get 
respectively\begin{eqnarray}(I \otimes R^{-1})(Z \otimes I) &+& 
(I \otimes R^{-1})(R^{-1} \otimes I)(I \otimes Q \otimes I)(I \otimes Z)
\nonumber\\
=(\tilde{R} \otimes I)(I \otimes Z)R^{-1} &+& (Z \otimes I)R^{-1} -
(I \otimes R^{-1})(R^{-1} \otimes I)(I \otimes Z)\end{eqnarray}and
\begin{eqnarray}
(\tilde{Z} \otimes I)\tilde{R} &-& (I \otimes Z)\tilde{R} - 
(I \otimes Q \otimes I)(I \otimes Z)\tilde{R}\nonumber\\=(I \otimes \tilde{R})
(\tilde{R} \otimes I)(I \otimes \tilde{Z}) &+& (I \otimes \tilde{R})
(I \otimes I \otimes Q)(Z \otimes I)\nonumber\\ &-& (\tilde{R} \otimes I)
(I \otimes \tilde{R})(Z \otimes I)
\end{eqnarray}
where we have used 
n and
\begin{eqnarray*}
\tilde{f}^{a}_{~b}(S(p^{n})a) = \tilde{f}^{c}_{~b}(S(p^{n}))\tilde{f}^{a}_{~c}(a) + \tilde{\eta}_{b}(S(p^{n}))\tilde{\eta}^{a}(a)
\end{eqnarray*}
and
\begin{eqnarray*}
\tilde{f}^{a}_{~b}(aS(p^{n})) = \tilde{f}^{c}_{~b}(a)\tilde{f}^{a}_{~c}(S(p^{n}))-\tilde{f}^{a}_{~b}(a_{(1)})\tilde{\eta}^{d}(a_{(2)})\tilde{\eta}_{d}(S(p^{n})),
\end{eqnarray*}
obtained from (26) for any $a \in {\cal A}$, to have the second equation. 
Multiplying both sides of (71) from the left by $(I \otimes R^{-1})$ on from 
the right by $R^{-1}$ and using (59)($\tilde{R} = R^{-1} + R^{-1}(I \otimes Q)$), we retrieve the equation (73). To investigate the equation (74), we have to 
multiply both sides of (72) from the left by $-(R \otimes I)$, to get
\begin{eqnarray*}
(\tilde{Z} \otimes I)\tilde{R} - (I \otimes Z)\tilde{R} &-& (I \otimes Q 
\otimes I)(I \otimes Z)\tilde{R}\nonumber\\
= - (R \otimes I)(I \otimes \tilde{R})(Z \otimes I) &-& (R \otimes I)
(I \otimes \tilde{R})(\tilde{R} \otimes I)(I \otimes Z).
\end{eqnarray*}
where we have used $\tilde{Z} = -RZ$ and $R\tilde{R}=I \otimes I + I \otimes Q$. Replacing $R$ by $\tilde{R} - I \otimes Q$ into the right hand side and 
using (67) and again $\tilde{Z} = -RZ$, we retrieve (74). Therefore, (71) 
implies (73) and (74).\\Finally the action of $\tilde{\eta}^{a} \circ S$ on both sides (52) gives\begin{eqnarray}(R \otimes I &-& I \otimes I \otimes I)
((I \otimes Z)Z - (Z \otimes I)Z)\nonumber\\
+ T \otimes I &-& 
(I \otimes \tilde{R})(\tilde{R} \otimes I)(I \otimes T)=0
\end{eqnarray}
where we have used the relations
\begin{eqnarray*}
\tilde{\eta}^{a}(S(ap^{n}))=\tilde{\eta}^{b}(S(p^{n}))\tilde{f}^{a}_{~b}
(S(a))
\end{eqnarray*}
and
\begin{eqnarray*}
\tilde{\eta}^{a}(S(p^{n}a))&=&\varepsilon(S(a))\tilde{\eta}^{a}(S(p^{n})) + 
\tilde{\eta}^{b}(S(a))\tilde{f}^{a}_{~b}(S(p^{n}))\\&-&\tilde{\eta}^{a}
(S(a_{(2)}))\tilde{\eta}^{b}(S(a_{(1)}))\tilde{\eta}_{b}(S(p^{n}))
\end{eqnarray*}
obtained from (28) for any $a \in {\cal A}$.\\We consider now the translations. From
\begin{eqnarray*}
\tilde{\eta}_{a}(S(p^{m})S(p^{n}))=\tilde{f}^{b}_{~a}(S(p^{m}))
\tilde{\eta}_{b}(S(p^{n})) - \tilde{\eta}_{a}(S(p^{k}))\tilde{\eta}^{b}
(S(\Lambda^{m}_{~k}))\tilde{\eta}_{b}(S(p^{n})),
\end{eqnarray*}
obtained from (27), we can see that the action of $\tilde{\eta}_{a} \circ S$ 
on both sides of (64) gives
\begin{eqnarray*}
(R - I \otimes I)((Q \otimes I)
\tilde{Z} - (Q \otimes I)ZQ -ZQ)=0 .
\end{eqnarray*}
From (65) and $\tilde{Z} =-RZ$ we can rewrite the left hand side of this 
equation as $-(R -I \otimes I)(R + I \otimes I + I \otimes Q)ZQ$ which vanishes by virtue of (62). Applying $f^{a}_{~b} \circ S$, $\tilde{f}^{a}_{~b} 
\circ S$ and $\tilde{\eta}^{a} \circ S$ on both sides of (64), we get 
respectively\begin{eqnarray}(R \otimes I &-& I \otimes I \otimes I)
((I \otimes Z)Z - (Z \otimes I)Z\nonumber\\&+& (I \otimes R^{-1})
(R^{-1} \otimes I)(I \otimes Q \otimes I)(I \otimes T))\nonumber\\+ 
(T \otimes I) &-& (I \otimes R^{-1})(R^{-1}
\otimes I)(I \otimes T)=0 ,\\(I \otimes R &-& I \otimes I \otimes I)((\tilde{Z} \otimes I)\tilde{Z} + (T \otimes I)Q -(I \otimes Q \otimes I)(I \otimes Z)
\tilde{Z}\nonumber\\&-&(I \otimes Q \otimes I)(I \otimes T) - (I \otimes Z)
\tilde{Z})\nonumber\\+(I \otimes T) &-& (\tilde{R} \otimes I)(I \otimes 
\tilde{R})(T \otimes I) = 0\end{eqnarray}and\begin{eqnarray}(I \otimes R) - 
I \otimes I \otimes I)((\tilde{Z} \otimes I)T - (I \otimes \tilde{Z})T)
\nonumber\\
- (Z \otimes I)T - (\tilde{R} \otimes I)(I \otimes Z)T=0
\end{eqnarray}
where we have used (25) to have (76),
\begin{eqnarray*}
\tilde{f}^{a}_{~b}(S(p^{m})S(p^{n})) = \tilde{f}^{c}_{~b}(S(p^{m}))
\tilde{f}^{a}_{~c}(S(p^{n})) + \tilde{\eta}_{b}(S(p^{m}))\tilde{\eta}^{a}
(S(p^{n}))\nonumber\\- \tilde{f}^{a}_{~b}(S(p^{k}))\tilde{\eta}^{c}
(S(\Lambda^{m}_{~k}))\tilde{\eta}_{c}(S(p^{n})) - \delta^{a}_{b}
\tilde{\eta}^{c}(S(p^{m}))\tilde{\eta}_{c}(S(p^{n}))
\end{eqnarray*}
obtained from (26) to have (77) and
\begin{eqnarray*}
\tilde{\eta}^{a}(S(p^{m})
(p^{m}))\tilde{f}^{a}_{~b}(S(p^{n})) - \tilde{\eta}^{a}(S(p^{k}))
\tilde{\eta}^{b}(S(\Lambda^{m}_{~k}))\tilde{\eta}_{~b}(S(p^{n}))
\end{eqnarray*}
obtained from (28) to have (78). The consistency of (75-78) requieres\\
{\bf Lemma (4,2)}: The relations (75-78) are consistent if
\begin{eqnarray}
T=-\tilde{R}T  &,\lambda \not=0,-2
\end{eqnarray}
$Proof.$ Multiplying both sides of (75) and (76) from the left by 
$(\tilde{R} \otimes I + I \otimes I \otimes I)$ and both sides of (77) and (78) from the left by $(I \otimes \tilde{R} + I \otimes I \otimes I)$ and using the 
relation $(\tilde{R} + I \otimes I)(R - I \otimes I) = 0$, we obtain the 
following consistency conditions \begin{eqnarray}(\tilde{R} \otimes I + 
I\otimes I \otimes I)(T \otimes I - (I \otimes \tilde{R})(\tilde{R} \otimes I)
(I \otimes T)) = 0,\\(\tilde{R} \otimes I + I\otimes I \otimes I)
(T \otimes I - (I \otimes R^{-1})(R^{-1} \otimes I)(I \otimes T)) = 0,\\
(I \otimes \tilde{R} + I\otimes I \otimes I)(I \otimes T - (\tilde{R} \otimes I)(I \otimes \tilde{R})(T \otimes I)) = 0
\end{eqnarray}
and
\begin{eqnarray}
(I \otimes \tilde{R} + I\otimes I \otimes I)((Z \otimes I)T + (\tilde{R} 
\otimes I)(I \otimes Z)T) = 0.
\end{eqnarray}
From (80) and (81), it follows
\begin{eqnarray}
(\tilde{R} \otimes I + I\otimes I \otimes I)((I \otimes \tilde{R})(\tilde{R} 
\otimes I)(I \otimes T)) \nonumber\\- (I \otimes R^{-1})(R^{-1} \otimes I)
(I \otimes T)) = 0,\end{eqnarray}and from (67) and $(I \otimes R^{-1})(R^{-1} 
\otimes I)(I \otimes \tilde{R}) = (\tilde{R} \otimes I)(I \otimes R^{-1})
(R^{-1} \otimes I)$ obtained by multiplying both sides of (68) from the left and from the right by$(I \otimes R^{-1})(R^{-1} \otimes I)$, we can rewrite (84) as\begin{eqnarray*}
((I \otimes \tilde{R})(\tilde{R} \otimes I) - 
(I \otimes R^{-1})(R^{-1} \otimes I))(I \otimes (\tilde{R}T + T))=0.
\end{eqnarray*}
Multiplying it by $(R \otimes I)(I \otimes R)$ and using $R\tilde{R} = 
I \otimes I + I \otimes Q$, we obtain $(I \otimes I \otimes T + T))=0$ 
which is equivalent, by virtue of (65), to $\lambda(2 + \lambda)
(\tilde{R} + I\otimes I)T = 0$ leading to (79) for $\lambda \not= 0$ or 
$-2$.From (79) and (67) we may see that the equation (82) is satisfied and it is also easy to see that the equation (83) is obtained by multiplying both sides 
of (72) from the right by T and by using (79).   Q.E.D.\\$Remark (4,1)$: 
note that for $\lambda=-2$ the characteristic equation (62) reduces to 
$(R - I \otimes I)^{2} = 0$ showing that $R$ does not possess negative  
eigenvalues and therefore, in this case, we can not construct from it 
 antisymmetric products. (In the following we shall assume the relation (79) 
for $\lambda = -2$).\\$Remark (4.2)$: For the case $\lambda =0$ we have 
$\tilde{R}=R=R^{-1}$ from which we see that the relations (80) and (81) are 
identical and can not imply (79). But by multiplying both sides of (64) by 
$\tilde{R} + I \otimes I$ and by using (62) we can deduce in this case the
following consi
stency condition
\begin{eqnarray}
C^{nm} - \Lambda^{n}_{~k} \Lambda^{m}_{~\ell}C^{k\ell}=0
\end{eqnarray}
where we have used (47) for $a=\Lambda$ and $C^{nm} = \tilde{R}^{nm}_{~k\ell} + \delta^{n}_{k} \delta^{m}_{\ell}$. In fact by applying $f^{a}_{~b} \circ S$ and $\tilde{f}^{a}_{~b}$ on both sides of (85) we get respectively
\begin{eqnarray*}
(\tilde{R} \otimes I + I \otimes I \otimes I)(T \otimes I) - 
(I \otimes R^{-1})(R^{-1} \otimes I)(I \otimes \tilde{R} + 
I \otimes I \otimes I)(I \otimes T)=0\end{eqnarray*} and
\begin{eqnarray*}
(I \otimes \tilde{R} + I \otimes I \otimes I)(I \otimes T) - (\tilde{R} 
\otimes I)(I \otimes \tilde{R})(\tilde{R} \otimes I +I \otimes I \otimes I)
(T \otimes I) = 0
\end{eqnarray*}
which are identical to the relations (81) and (82) by virtue of (67). To get the relation (83), we apply $\tilde{\eta}^{a} \circ S$ on both sides of (85) to 
have
\begin{eqnarray*}
(Z \otimes I)(\tilde{R} + I \otimes I)T + (\tilde{R} \otimes I)(I \otimes Z)
(\tilde{R} + I \otimes I)T=0
\end{eqnarray*}
and (36). By virtue of (72) we can see that this relation is equivalent to (83). Therefore for $\lambda = 0$ the consistency condition is not $RT = -T$ as 
assumed in remark (3.11) of Ref.[9] but (85).\\
Comparing now (75) with (76), we obtain
\begin{eqnarray}
(R \otimes I &-& I \otimes I \otimes I)(I \otimes R^{-1})(R^{-1} \otimes I)
(I \otimes Q \otimes I)(I \otimes T)\nonumber\\
 &-& (I \otimes R^{-1})(R^{-1} \otimes I)(I \otimes T)= -(I \otimes \tilde{R})
(\tilde{R} \otimes I)(I \otimes T).
\end{eqnarray}
Replacing $R$ by $\tilde{R} - I \otimes Q$ and using $(\tilde{R} \otimes I)
(I \otimes R^{-1})(R^{-1} \otimes I)=(I \otimes R^{-1})(R^{-1} \otimes I)
(I \otimes \tilde{R})$, (65) and (79), we can write the left hand side as 
\begin{eqnarray*}
-(I \otimes R^{-1})(R^{-1} \otimes I)(I \otimes I \otimes I + 
I \otimes I \otimes Q + I \otimes Q \otimes I + I \otimes I \otimes Q^{2})
(I \otimes T).
\end{eqnarray*}
which can be identified to the right hand side of (86) by using (
59) and (65). Therefore, (75) is equivalent to (76).\\
$Remark (4.3)$: For $\lambda =0$ the relations (75) and (76) are identical but 
for $\lambda = -2$ we must assume the relation (79) to have the equivalency 
between (75) and (76).\\
To investigate the equation(77), we multiply both sides of (71) from the left 
by $(R \otimes I)$ and from the right by $Z$ to obtain
\begin{eqnarray*}
(\tilde{Z} \otimes I)\tilde{Z}&+&(I \otimes \tilde{R})(I \otimes \tilde{Z})
\tilde{Z}\\=(R\otimes I)(I \otimes R)(Z \otimes I)Z&+&(R \otimes I)
(I \otimes R)(\tilde{R} \otimes I)(I \otimes Z)Z
\end{eqnarray*}
where we have used $\tilde{Z} = -RZ$, $R\tilde{R} = I \otimes I + I \otimes Q = \tilde{R}R = I \otimes I + Q \otimes I$ (due to (62) and (65)). Now, 
multiplying both sides from the left by $I \otimes R - I \otimes I \otimes I$, 
we obtain
\begin{eqnarray*}
(I \otimes R - I \otimes I \otimes I)((\tilde{Z} \otimes I)\tilde{Z} - 
(I \otimes \tilde{Z})\tilde{Z}) = \\(I \otimes R - I \otimes I \otimes I)
(R \otimes I)(I \otimes \tilde{R} \otimes I)(I \otimes Z)Z)
\end{eqnarray*}
where we have used (68) and $(R - I \otimes I)\tilde{R} = -(R - I \otimes I)$.
Taking into account this equation, (77) can be written as
\begin{eqnarray*}
(I \otimes R &-& I \otimes I \otimes I)((R \otimes I)(I \otimes R)((Z \otimes I)Z + (\tilde{R} \otimes I)(I \otimes Z)Z)\\ &+& (I \otimes I \otimes Q)
(T \otimes I)-(I \otimes Q \otimes I)(I \otimes T))\\&+& I \otimes T - 
(\tilde{R} \otimes I)(I \otimes \tilde{R})(T \otimes I) =0.
\end{eqnarray*}
Multiplying both sides from the left by $-(I \otimes R^{-1})(R^{-1} \otimes I)$ and using $(I \otimes R^{-1})(R^{-1} \otimes I)(I \otimes R)=(R \otimes I)
(I \otimes R^{-1})(R^{-1} \otimes I)$, we obtain
\begin{eqnarray*}
(R \otimes I &-& I \otimes I \otimes I)((I \otimes Z)Z -(Z \otimes I)Z\\ &-&
(I \otimes R^{-1})(R^{-1} \otimes I)(I \otimes I \otimes Q)(T \otimes I)\\&+&
(I \otimes R^{-1})(R^{-1} \otimes I)(I \otimes Q \otimes I)(I \otimes T))\\&-&
(I \otimes R^{-1})(R^{-1} \otimes 
I)(I \otimes T)\\&+&(I \otimes R^{-1})(R^{-1} \otimes I)(\tilde{R} \otimes I)
(I \otimes \tilde{R})(T \otimes I)=0=\\-(R \otimes I &-& I \otimes I \otimes I)
(I \otimes R^{-1})(R^{-1} \otimes I)(I \otimes I \otimes Q)(T \otimes I)\\ - 
T \otimes I&+&(I \otimes R^{-1})(R^{-1} \otimes I)(\tilde{R} \otimes I)
(I \otimes \tilde{R})(T \otimes I)\end{eqnarray*}where we have used (76) 
to have the third hand side which gives, by using ( 56) into the third factor 
of the last term,\begin{eqnarray*}&-&(R \otimes I - I \otimes I \otimes I)
(I \otimes R^{-1})(R^{-1} \otimes I)(I \otimes I \otimes Q)(T \otimes I)\\ &+& 
(I \otimes R^{-1})(I \otimes I \otimes Q)(T \otimes I)\\&+& (I \otimes R^{-1})
(R^{-1} \otimes I)(I \otimes \tilde{R})(I \otimes I \otimes Q)(T \otimes I) = 0.\end{eqnarray*}For $\lambda = 0$ this relation is trivial and for $\lambda 
\not= 0$ we use (79) and (59) into the second term and (56) into the third term to obtain 
\begin{eqnarray*}
&-&(R \otimes I - I \otimes I \otimes I)(I \otimes R^{-1}\otimes I \otimes Q)
(T \otimes I)\\ &-& (I \otimes R^{-1})(R^{-1} \otimes I)(I \otimes I \otimes Q)
(T \otimes I)\\&+& (I \otimes R^{-1})(R^{-1} \otimes I)(I \otimes R)
(I \otimes I \otimes Q)(T \otimes I) = 0
\end{eqnarray*}
which is satisfied by virtue of $(I \otimes R^{-1})(R^{-1} \otimes I)
(I \otimes R) = (R \otimes I)(I \otimes R^{-1})(R^{-1} \otimes I)$. 
Therefore, (77) is also equivalent to (75) and (76).\\
Now, we consider the consistency of the commutation rules (64) and the braiding properties of the product of three generators of ${\cal B}$. After the insertion of (64) into the first term of the right hand side of (64) it self, we get
\begin{eqnarray*}
p \top p = R^{2}(p \top p) - (R + I \otimes I)(R - I \otimes I)Zp + \\
(R + I \otimes I)T - (\Lambda \top \Lambda)(R + I \otimes I)T = \\
p \top p - (I \otimes Q)((R - I \otimes I)((p \top p)-Zp) + T - 
(\Lambda \top \Lambda)T)
\end{eqnarray*}
where we have used into the second hand side (42) for $a=\Lambda$ to have the
 last term and (62), (56) and (79) to have the third hand side. We observe from (64) that the second term of the third hand side vanishes implying the consistency of these commutation rules. For the product of three translations, we have
\begin{center}
\begin{eqnarray*}
p \top p \top p = (p \top p)\top p =(R \otimes I)(p \top (p \top p)) + 
r \top p = \\(R \otimes I)((I \otimes R)((p \top p) \top p) + p \top r) + 
r \top p =\\(R \otimes I)(I \otimes R)(R \otimes I)(p \top p \top p) + 
(R \otimes I)(I \otimes R)(r \top p) +\\ (R \otimes I)(p \top r) + 
r \top p
\end{eqnarray*}
\end{center}
where $r = -(R - I \otimes I)Zp + T - (\Lambda \top \Lambda)T$. On the other 
hand
\begin{center}
\begin{eqnarray*}
p \top p \top p = p \top (p \top p) = (I \otimes R)((p \top p) \top p) + 
p \top r= \\(I \otimes R)((R \otimes I)(p \top (p \top p)) + r \top p) + 
p \top r =\\(I \otimes R)(R \otimes I)(I \otimes R)(p \top p \top p) + 
(I \otimes R)(R \otimes I)(p \top r) + \\(I \otimes R)(r \top p) + p \top r.
\end{eqnarray*}
\end{center}
g-Baxter equation (66) implies \begin{eqnarray}A(r \top p) = B(p \top r)
\end{eqnarray}
where\begin{eqnarray*}A=(R \otimes I)(I \otimes R) - I \otimes R + I \otimes 
I \otimes I\\B=(I \otimes R)(R \otimes I) - R \otimes I + I \otimes I \otimes I.\end{eqnarray*}Replacing into both sides of (87) $r$ by $-(R-I \otimes I)Zp 
+ T - (\Lambda \top \Lambda)T$, we obtain\begin{eqnarray}-A(R \otimes I - 
I \otimes I \otimes I)(Z \otimes I)(p \top p) + A(T \otimes I)p -A
(\Lambda \top \Lambda \top p)T=\nonumber\\-B(I \otimes R - I \otimes 
I \otimes I)(I \otimes Z)(p \top p) +B(I \otimes T)p -B(p \top 
\Lambda \top \Lambda)T.\end{eqnarray}Using (52)$(p \top \Lambda = \tilde{R}
(\Lambda \top p) + Z\Lambda - (\Lambda \top \Lambda)Z)$, we can rewrite the 
third term of the right hand side of (88) as\begin{eqnarray}B(p \top \Lambda 
\top \Lambda)T &=& B(\tilde{R} \otimes I)(\Lambda \top p \top \Lambda)T
\nonumber\\ + B(Z \otimes I)(\Lambda \top \Lambda)T &-& B(\Lambda \top 
\Lambda \to
p \Lambda)T=\nonumber\\B(\tilde{R} \otimes I)(I \otimes \tilde{R})
(\Lambda \top \Lambda \top p)T &+&B(\tilde{R} \otimes I)(I \otimes Z)
(\Lambda \top \Lambda)T\nonumber\\ - B(\tilde{R} \otimes I)(\Lambda \top 
\Lambda \top \Lambda)(I \otimes Z)T &+& B(Z \otimes I)(\Lambda \top \Lambda)T - B(\Lambda \top \Lambda \top \Lambda)(Z \otimes I)T.\nonumber\\\end{eqnarray}A 
straightforward computation shows that\begin{eqnarray}B(\tilde{R} \otimes I)
(I \otimes \tilde{R}) = A + (I \otimes I \otimes Q)(\tilde{R} \otimes I + 
I \otimes I \otimes I)\end{eqnarray}where we have used $R\tilde{R} = 
I \otimes I + I \otimes Q$ and (56). Multiplying both sides of (90) from 
the right by  $(\Lambda \top \Lambda \top p)T$ and using (47) for $a=\Lambda$
 and (79), we get\begin{eqnarray*} A(\Lambda \top \Lambda \top p)T = 
B(\tilde{R} \otimes I)(I \otimes \tilde{R})(\Lambda \top \Lambda \top p)T.
\end{eqnarray*}
We see from (90) that for $\lambda = 0$ this relation is satisfied without 
using (79). On the other hand by using
$\Lambda$, we can rewrite the third and the fifth term of third hand of (89) as \begin{eqnarray*}(\Lambda \top \Lambda \top \Lambda)B(Z \otimes I + 
(\tilde{R}\otimes I)(I \otimes Z))T.\end{eqnarray*}Using $\tilde{Z} = -RZ$, 
$R\tilde{R} = \tilde{R} R = I \otimes I + I \otimes Q$ and (56) we get 
\begin{eqnarray}
B(Z \otimes I + (\tilde{R} \otimes I)(I \otimes Z)T =\nonumber\\
-(I \otimes R - I \otimes I \otimes I)(\tilde{Z} \otimes I - I \otimes 
\tilde{Z})T + (Z \otimes I)T\nonumber\\
 + (\tilde{R} \otimes I)(I \otimes Z)T 
\end{eqnarray}
which vanishes by virtue of (78). Therefore (88) reduces to 
\begin{eqnarray}
-A(R \otimes I &-& I \otimes I \otimes I)(Z \otimes I)(p \top p) + A(T \otimes I)p =\nonumber\\
 -B(I \otimes R &-& I \otimes I \otimes I)(I \otimes Z)(p \top p) + B(I \otimes T)p\nonumber\\
- B(Z \otimes I &+& (\tilde{R} \otimes I)(I \otimes Z))(\Lambda \top \Lambda)T.
\end{eqnarray}
Using now \begin{eqnarray*}A(R \otimes I - I \otimes I \otimes I)=B(I \otimes R - I \otimes I \otimes I)= \Lambda_{3},
\end{eqnarray*}
we may rewrite (92) as
\begin{eqnarray}
\Lambda_{3}(Z \otimes I - I \otimes Z)(p \top p) - A(T \otimes I)p + B(I \otimes T)p - \nonumber\\
B(Z \otimes I + (\tilde{R} \otimes I)(I \otimes Z))(\Lambda \top \Lambda)T=0.
\end{eqnarray}
From the properties of $\Lambda_{3}$ and $(R - I \otimes I)\tilde{R} =
 -(R - I \otimes I)$, we have
\begin{eqnarray*}
\Lambda_{3}(Z \otimes I - I \otimes Z)(p \top p) = \Lambda_{3}(Z \otimes I + 
(\tilde{R} \otimes I)(I \otimes Z))(p \top p) =\\B(I \otimes R)(Z \otimes I +
 (\tilde{R} \otimes I)(I \otimes Z))(p \top p)-B(Z \otimes I + 
(\tilde{R} \otimes I)(I \otimes Z))(p \top p)
\end{eqnarray*}
Using (64) into the second term of the third hand side, we obtain
\begin{center}\begin{eqnarray*}
B(I \otimes R)(Z \otimes I + (\tilde{R} \otimes I)(I \otimes Z))
(p \top p)-B(Z \otimes I + (\tilde{R} \otimes I)(I \otimes Z))R(p \top p)\\
+B(Z \otimes I + (\tilde{R} \otimes I)(I \otimes Z))(R - I \otimes I)Zp -
B(Z \otimes Z))T\\+B(Z \otimes I + (\tilde{R} \otimes I)(I \otimes Z))
(\Lambda \top \Lambda)T.\end{eqnarray*}\end{center}The first and the second 
term vanish by virtue of (71) and the fourth term vanishes by virtue of (91). 
Therefore (93) reduces to \begin{eqnarray*}+B(Z \otimes I + 
(\tilde{R} \otimes I)(I \otimes Z))(R - I \otimes I)Zp - A(T \otimes I)p + 
B(I \otimes T)p\end{eqnarray*}which can be written, after the make of use of 
(71), as\begin{center}\begin{eqnarray*}+B(I \otimes R - I \otimes I \otimes I)
(Z \otimes I + (\tilde{R} \otimes I)(I \otimes Z))Zp - A(T \otimes I)p + 
B(I \otimes T)p =\\\Lambda_{3}(Z \otimes I + (\tilde{R} \otimes I)
(I \otimes Z))Zp - A(T \otimes I)p + B(I \otimes T)p =\\
A(R \otimes I - I \otimes I \otimes I)((Z \otimes I)Z - (I \otimes Z)Z)p -
 A(T \otimes I)p + B(I \otimes T)p\end{eqnarray*}\end{center}which reduces, 
due to (75), to \begin{eqnarray*}-A(I \otimes \tilde{R})(\tilde{R} \otimes I)
(I \otimes T)p + B(I \otimes T)p
\end{eqnarray*}
which vanishes according to the same calculation of (90).\\
The same procedure is applied to $p\top p \top \Lambda$ to have
\begin{center}
\begin{eqnarray*}
p\top p \top \Lambda =p \top (p \top \Lambda) = (I \otimes \tilde{R})
((p \top \Lambda) \top p) + (I \otimes Z)(p \top \Lambda) -\\
 ((p \top \Lambda) \top \Lambda )Z=\\(I \otimes \tilde{R})(\tilde{R} \otimes I)
(\Lambda \top (p \top p)) +(I \otimes \tilde{R})(Z \otimes I)(\Lambda \top p)
 - (I \otimes \tilde{R})(\Lambda \top \Lambda \top p)Z +\\(I \otimes Z)\
tilde{R}(\Lambda \top p) + (I \otimes Z)Z\Lambda - (I \otimes Z)
(\Lambda \top \Lambda)Z -(\tilde{R} \otimes I)(\Lambda \top  
(p \top \Lambda))Z - \\ (Z \otimes I)(\Lambda \top \Lambda)Z + 
(\Lambda \top \Lambda \top \Lambda)(Z \otimes I)Z =\\  (I \otimes \tilde{R})
(\tilde{R} \otimes I)(I \otimes R)(\Lambda \top p \top p) - \\(I \otimes 
\tilde{R})(\tilde{R} \otimes I)(I \otimes R - I \otimes I \otimes I)
(I \otimes Z)(\Lambda \top p) +\\ (I \otimes \tilde{R})(\tilde{R} \otimes I)
(I \otimes\tilde{R})(\tilde{R} \otimes I)(\Lambda \top \Lambda \top \Lambda)
(I \otimes T) +\\ (I \otimes \tilde{R})(Z \otimes I)(\Lambda \top p) -
(I \otimes \tilde{R})(\Lambda \top \Lambda \top p)Z + (I \otimes Z)\tilde{R}
(\Lambda \top p) + \\ (I \otimes Z)Z\Lambda -(I \otimes Z)
(\Lambda \top \Lambda)Z - (\tilde{R} \otimes I)(I \otimes \tilde{R})
(\Lambda \top \Lambda \top p)Z -\\ (\tilde{R} \otimes I)(I \otimes Z)
(\Lambda \top \Lambda)Z + \\(\tilde{R} \otimes I)(\Lambda \top \Lambda \top 
\Lambda)(I \otimes Z)Z - (Z \otimes I)(\Lambda \top \Lambda)Z + 
(\Lambda \top \Lambda \top \Lambda)(Z \otimes I)Z.
\end{eqnarray*}
\end{center}
on the other hand,
\begin{center}
\begin{eqnarray*}
p \top p \top \Lambda = (p \top p) \top \Lambda =(R \otimes I)
(p \top (p \top \Lambda)) -\\ (R \otimes I - I \otimes I \otimes I)
(Z \otimes I)(p \top \Lambda) + (T \otimes I)\Lambda - 
(\Lambda \top \Lambda \top \Lambda)(T \otimes I) =\\(R \otimes I)
(I \otimes \tilde{R})((p \top \Lambda) \top p) + (R \otimes I)
(I \otimes Z)(P \top \Lambda) - (R \otimes I)((p \top \Lambda) \top \Lambda)Z 
-\\(R \otimes I -I \otimes I \otimes I)(Z \otimes I) \tilde{R}(\Lambda \top p) -(R \otimes I -I \otimes I \otimes I)(Z \otimes I)Z\Lambda +\\
(R \otimes I -I \otimes I \otimes I)(Z \otimes I)(\Lambda \top \Lambda)Z +
(T \otimes I)\Lambda - (\Lambda \top \Lambda \top \Lambda)(T \otimes I) =\\
(R \otimes I)(I \otimes \tilde{R})(\tilde{R} \otimes I)
(\Lambda \top p \top p) +(R \otimes I)(I \otimes \tilde{R})(Z \otimes I)
(\Lambda \top p) -\\(R \otimes I)(I \otimes \tilde{R})
(\Lambda \top \Lambda \top p)Z + (R \otimes I)(I \otimes Z)\tilde{R}
(\Lambda \top p) +\\(R \otimes I)(I \otimes Z)Z\Lambda - 
(R \otimes I)(I \otimes Z)(\Lambda \top \Lambda)Z - (R \otimes I)(\tilde{R} 
\otimes I)(\Lambda \top p \top \Lambda)Z -\\ (R \otimes I)(Z \otimes I)
(\Lambda \top \Lambda)Z + (R \otimes I)(\Lambda \top \Lambda \top \Lambda)(Z \otimes I)Z -\\(R \otimes I -I \otimes I \otimes I)(Z \otimes I) \tilde{R}
(\Lambda \top p) -(R \otimes I \otimes I)Z\Lambda +\\
(R \otimes I -I \otimes I \otimes I)(Z \otimes I)(\Lambda \top \Lambda)Z +
(T \otimes I)\Lambda - (\Lambda \top \Lambda \top \Lambda)(T \otimes I). \\
\end{eqnarray*}
\end{center}
Using again (52), we can replace $(R \otimes I)(\tilde{R} \otimes I)
(\Lambda \top p \top \Lambda)Z$ by:
\begin{eqnarray*}
 (R \otimes I)(\tilde{R} \otimes I)(I \otimes \tilde{R})
(\Lambda \top \Lambda \top p)Z + (R \otimes I)(\tilde{R} \otimes I)
(I \otimes Z)(\Lambda \top \Lambda)Z\\ - (R \otimes I)(\tilde{R} \otimes I)
(\Lambda \top \Lambda \top \Lambda)(I \otimes Z)Z.
\end{eqnarray*}
We remark that (69) implies the egality between the coefficients multiplying 
$\Lambda \top p \top p$.\\ 
For the coefficient multiplying $\Lambda \top p$, we must have the following 
relation\begin{eqnarray*}(R \otimes I)(I \otimes \tilde{R})(Z \otimes I) + 
(R \otimes I)(I \otimes Z)\tilde{R} - (R \otimes I -I \otimes I \otimes I)
(Z \otimes I) \tilde{R} =\\-(I \otimes \tilde{R})(\tilde{R} \otimes I)
(I \otimes R - I \otimes I \otimes I)(I \otimes Z) + (I \otimes \tilde{R})
(Z \otimes I) + (I \otimes Z)\tilde{R}
\end{eqnarray*}which can be written under the form
\begin{eqnarray*}
(R \otimes I -I \otimes I \otimes I)((I \otimes \tilde{R})(Z \otimes I) - 
(Z \otimes I)\tilde{R} \\+ (I \otimes Z)\tilde{R} + (I \otimes \tilde{R})
(\tilde{R} \otimes I)(I \otimes Z)) =0
\end{eqnarray*}
where we have used (69). By using $(R - I \otimes I)\tilde{R} =
 -(R - I \otimes I)$ we can replace the third term of the second factor by 
$-(\tilde{R} \otimes I)(I \otimes Z)\tilde{R}$ to see that this relation is 
satisfied by virtue of (72). It is also easy to see that the coefficients 
multiplying $\Lambda$ are precisely the relation (75) and those multiplying 
$(\Lambda \top \Lambda)Z$ and $(\Lambda \top \Lambda \top p)Z$ are equal in 
virtue of $R\tilde{R} = I \otimes I + I \otimes Q$ and (56). Using (42) and 
(47) for $a = \Lambda$, we obtain for the coefficients multiplying 
$\Lambda \top \Lambda \top \Lambda$ the same consistency of the braiding of 
$p \top p \top \Lambda$.\\  
Finally by using $\Lambda \top \Lambda = R(\Lambda \top \Lambda)R^{-1}$, 
obtained from (42) for $a = \Lambda$, and (71) we can show by a similar 
way the consistency of the braiding of $p \top \Lambda \top \Lambda$.\\
 Then from the results of this section, it follows\\ {\bf Theorem (4,1)}:Let  
${\cal G}$ be an inhomogeneous quantum group, as defined in section 2, with the following commutation rules
\begin{eqnarray*}
R(\Lambda \top \Lambda)&=&(\Lambda \top \Lambda)R,\\p \top \Lambda &=& 
\tilde{R}(\Lambda \top p) + Z\Lambda - (\Lambda \top \Lambda)Z,\\p \top p 
&=& R(p \top p) - (R - I \otimes I)Zp + T - (\Lambda \top \Lambda)T,
\end{eqnarray*}
then
\begin{eqnarray*}
\tilde{R} = R + I \otimes Q &=& R^{-1} + R^{-1}(I \otimes Q),\\Q=
\lambda I~~~\lambda \in {\cal C}&,&~~~\lambda \not= -1,\\T=
-\tilde{R}T~~~&for& \lambda \not= 0\\or~~ (\tilde{R} + I \otimes I)T - 
(\Lambda \top \Lambda)(\tilde{R} + I \otimes I)T = 0~~
&for& \lambda =0,\\ (R \otimes I)(I \otimes R)(R \otimes I)&=&(I \otimes R)
(R \otimes I)(I \otimes R),\\(Z \otimes I)R &+& (\tilde{R} \otimes I)
(I \otimes Z)R\\=(I \otimes R)(Z \otimes I) &+& (I \otimes R)
(\tilde{R} \otimes I)(I \otimes Z), \\(R \otimes I - I \otimes I \otimes I)&&
((I \otimes Z)Z - (Z \otimes I)Z)\\+ T \otimes I &-& (I \otimes \tilde{R})
(\tilde{R} \otimes I)(I \otimes T)=0\end{eqnarray*}and\begin{eqnarray*}
(I \otimes R) - I \otimes I \otimes I)((\tilde{Z} \otimes I)T - 
(I \otimes \tilde{Z})T)- (Z \otimes I)T - (\tilde{R} \otimes I)(I \otimes Z)T=0.\end{eqnarray*}with\begin{eqnarray*}\tilde{Z} = -RZ\end{eqnarray*}
\section{\bf Dicussions and Conclusions}We end this paper by noticing that:\\1. Let us recall that (10) and (11) imply that $f^{n}_{\ell}$, $f^{nm}_{k\ell}$, 
$f^{nm}_{~\ell}$ and $f^{n}_{k\ell}$ are linear functionals on ${\cal B}$ but not $\tilde{f}^{m}_{~k}, \tilde{\eta}^{m}$ and $\tilde{\eta}_{k}$ . 
Then their action on both sides of inhomogeneous commutation
not justified but a tedious and straightforward computation shows that they 
give the same results that those obtained directly by using $f^{nm}_{k\ell}$, 
$f^{nm}_{~\ell}$ and $f^{n}_{k\ell}$ .\\  2. From (35), (36) and (24), 
we see that ${\cal A} \ni a \rightarrow \tilde{\rho}(a)= \left(\begin{array}{cc}                                               \tilde{f}(a)  &  \tilde{\eta}(a)\\                                               0   &  \varepsilon(a)                                                \end{array}                                               \right) \in M_{N}(C)$ is a unital antihomomorphism.\\3. Replacing $a$ and $b \in {\cal A}$ respectively by $S(a)$ and $S(b)$ into (36) and by 
setting $\eta^{k}= \tilde{\eta}^{k} \circ S$, we get\begin{eqnarray}
\eta^{n}(ab)= \eta^{n}(a) \varepsilon(b)+ \tilde{f}^{n}_{~m}(S(a)) \eta^{m}(b),  &a,b \in {\cal A}.\end{eqnarray}Then, although the formalism presented above is quite different from those of Ref.[9], we arrive 
to similar results for the commutation rules between the elements of ${\cal A}$ and the translations (48) and (94). These formulae become identical with those 
of Ref.[9] if $\tilde{f}^{n}_{~m}(S(a))=f^{n}_{~m}(a)$ for any $a \in {\cal A}$. The latter condition is satisfied for $\tilde{\eta}_{n}(p^{a})=0$ and it is 
true, as a consequence of (56) and (57), for any $a \in {\cal B}$. In this case (34), (35) and (94) are also true for any $a$ and $b \in {\cal B}$ and (94) 
can be combined with (34) to see that 
${\cal A} \ni a \rightarrow \rho(a)= \left(\begin{array}{cc}                                               f(a)  &  \eta(a) \\                                               0   &  \varepsilon(a)                                                \end{array}                                         \right) \in M_{N}(C)$ is aunital homomorphism. In this case the $R$-matrix is subject to a strong
constraint $R^{2}=I^{\otimes 2}$ as in Ref.[9].\\4. Using (85) and setting 
$T=-2T'$ we can re
   write, f
ion rules between the translations as\begin{eqnarray*}(R^{nm}_{k\ell} - 
\delta^{n}_{k}\delta^{m}_{\ell})(p^{k}p^{\ell}-Z^{kl}_{q}p^{q} + 
T'^{k\ell} -\Lambda^{k}_{~q}\Lambda^{\ell}_{~p}T'^{qp}) = 0
\end{eqnarray*}
which are identical to those of Ref.[9]. In this case  one can also see 
that (75) reduce to the condition of the existence of covariant differential 
calculus on a quantum Minkowski space found in Ref.[12].\\

{\bf Acknowledgments.} We are grateful to M. Dubois-Violette, A. B. Hammou and  A.E.K. Yanallah for valuable discussions\\

{\bf References.}\\
1- O. Ogiesvestsky, W. B. Schmidke, J. Wess and B. Zumino, Commun. Math. Phys. 150(1992)495.\\
2- S. Majid, J. Math. Phys. 34(1993)2045.\\
3- M. Chaichian and A.P. Demichev, Phys. Lett. B304(1993)220.\\
4- V. K. Dobrev, J. Phys. A: Math. Gen. 26(1993),1317.\\ 
5- J. Lukierski, A. Nowicki, H. Ruegg and V.N. Tolstoy, Phys. Lett. B264(1991)
331.\\ 
6- J. Lukierski and A. Nowicki, Phys. Lett. B279(1992)299, M. Schliecker, 
W. Weich and R. Weixler, Z. Phys. C. Particles and Fields 53(1992)79, Lett. 
Math. Phys. 27(1993)217, L. Castellani, Lett. Math. Phys. 30(1994)233.\\
7- S. Zakrzewski, J. Phys. A:Math. Gen. 27(1994),2075.\\
8- L. Castellani, Commun. Math. Phys. 171(1995)383, P. Aschieri and L. Castellani,  Int. J. Mod. Phys. A11(1996)4513.\\
9- P. Podles and S.L. Woronowicz, Commun. Math. Phys. 185(1997)325, Commun. Math. Phys. 178(1996)61  , "Inhomogeneous quantum groups" In Proceedings of First Caribbean School of Mathematics and Physics in Guadeloupe, 1993. M. Chaichian and A. Demichev, "Quantum Groups" World Scientific (1996)\\
10- S.L. Woronowicz, Commun. Math. Phys. 122(1989)125.\\
11- U. Carow-Watamura,M Schliecker, S. Watamura and W. Weich, Commun. Math. Phys.  142(1991)605. M. Lagraa, Int. J. Mod. Phys. A11(1996)699.\\
12- P. Podles, Commun. Math. Phys. 181(1996)569.

\end{document}